\documentclass[amssymb,11pt]{amsart}
\usepackage{amscd,amssymb,euscript,amsthm}
\theoremstyle{plain}
\newtheorem{thm}{Theorem}[section] 

\newtheorem{lem}[thm]{Lemma}
\theoremstyle{definition}
\newtheorem{defn}[thm]{Definition}
\theoremstyle{remark}
\newtheorem{rem}[thm]{Remark}

\numberwithin{equation}{section}
\def\<{\left<}
\def\>{\right>}
\begin{document}
\title[Existence of $E_0$-semigroups]{On the existence of $E_0$-semigroups}
\author{William Arveson}
%
%
\address{Department of Mathematics,
University of California, Berkeley, CA 94720}
\email{arveson@math.berkeley.edu}
\maketitle

\begin{abstract}
Product systems are the classifying structures for 
semigroups of endomorphisms of $\mathcal B(H)$, 
in that two $E_0$-semigroups are cocycle conjugate iff their 
product systems are isomorphic.  Thus it is important to know 
that every abstract product system is associated with an 
$E_0$-semigrouop.  This was first proved more than fifteen years ago  
by rather indirect methods.  
Recently, Skeide has given a more direct proof.  In this 
note we give yet another proof by an elementary 
construction.  
\end{abstract}
\maketitle

\section{Formulation of the result}\label{S:ss}
There were two proofs of the above fact 
\cite{arvIV}, \cite{lieb1} (also see \cite{arvMono}), both of 
which involved substantial analysis.  
In a recent paper, Michael Skeide \cite{skeid1} gave a more direct 
proof.  In this note we present an elementary method for constructing 
an essential representation of any product system.  Given the basic 
correspondence between $E_0$-semigroups and essential representations, 
the existence of an appropriate $E_0$-semigroup follows.  

Our terminology follows the monograph \cite{arvMono}.  
Let $E=\{E(t):t>0\}$ be a product system and choose 
a unit vector $e\in E(1)$. {\em $e$ will be fixed throughout.}  We consider 
the Fr\'echet space of all 
Borel - measurable sections $t\in(0,\infty)\mapsto f(t)\in E(t)$ that are 
locally square integrable 
\begin{equation}\label{ssEq1}
\int_0^T\|f(\lambda)\|^2\,d\lambda<\infty, \qquad T>0.  
\end{equation}

\begin{defn}
A locally $L^2$ section $f$ is said to be {\em stable} if there is 
a $\lambda_0>0$ such that for almost every $\lambda\geq\lambda_0$, one has 
$$
f(\lambda+1)=f(\lambda)\cdot e .  
$$
\end{defn}
  
Note that a stable section $f$ satisfies $f(\lambda+n)=f(\lambda)\cdot e^n$ a.e. for 
all $n\geq 1$ whenever $\lambda$ is sufficiently large.  
The set of all stable sections is a vector space $\mathcal S$, and for any two 
sections $f,g\in \mathcal S$, $\langle f(\lambda+n),g(\lambda+n)\rangle$ becomes independent 
of $n\in \mathbb N$ (a.e.) when $\lambda$ is sufficiently large.  Thus we can define 
a positive semidefinite inner product on $\mathcal S$ as follows 
\begin{equation}\label{ssEq2}
\langle f,g\rangle=\lim_{n\to\infty}\int _n^{n+1}\langle f(\lambda),g(\lambda)\rangle\,d\lambda
=\lim_{n\to\infty}\int_0^1\langle f(\lambda+n),g(\lambda+n)\rangle\,d\lambda.
\end{equation}
Let $\mathcal N$ be the subspace of $\mathcal S$ consisting  
of all sections $f$ that vanish eventually, in that 
for some $\lambda_0>0$ one has $f(\lambda)=0$ for almost all $\lambda\geq \lambda_0$.  
One finds that $\langle f,f\rangle=0$ iff $f\in \mathcal N$.  Hence $\langle\cdot,\cdot\rangle$ 
defines an inner product on the quotient $\mathcal S/\mathcal N$, and 
its completion 
becomes a Hilbert space $H$ with respect to the inner product (\ref{ssEq2}).  
Obviously, $H$ is separable.  

There is a natural representation of $E$ on $H$.  Fix $v\in E(t)$, 
$t>0$.  For every stable section $f\in\mathcal S$, let $\phi_0(v)f$ be 
the section 
$$
(\phi_0(v)f)(\lambda)=
\begin{cases}
v\cdot f(\lambda-t), \quad &\lambda>t, \\
0, &0<\lambda\leq t.
\end{cases}
$$
Clearly $\phi_0(v)\mathcal S\subseteq \mathcal S$.  Moreover, $\phi_0(v)$ maps 
null sections into null sections, hence it induces a linear operator 
$\phi(v)$ on $\mathcal S/\mathcal N$.    The mapping 
$(t,v),\xi\in E\times \mathcal S/\mathcal N\mapsto \phi(v)\xi\in H$ 
is obviously Borel-measurable, and it is easy to check that 
$\|\phi(v)\xi\|^2=\|v\|^2 \cdot\|\xi\|^2$ (see Section \ref{S:re}
for details).   
Thus we obtain a representation $\phi$ of $E$ on the completion $H$ of $\mathcal S/\mathcal N$ 
by closing the densely defined operators
$\phi(v)(f+\mathcal N)=\phi_0(v)f+\mathcal N$, $v\in E(t)$, $t>0$, $f\in \mathcal S$.

\begin{thm}\label{esThm1}  $\phi$ is an essential representation of $E$ 
on $H$.
\end{thm}  

By Proposition 2.4.9 of \cite{arvMono}, there is an $E$-semigroup 
$\alpha=\{\alpha_t: t\geq 0\}$  
that acts on $\mathcal B(H)$ and is associated with $\phi$ by way of  
\begin{equation}\label{esEq0}
\alpha_t(X)=\sum_{n=1}^\infty \phi(e_n(t))X\phi(e_n(t))^*,\qquad X\in\mathcal B(H),
\quad t>0, 
\end{equation}
$e_1(t), e_2(t),\dots$ denoting an arbitrary orthonormal basis for $E(t)$.   
Since $\phi$ is essential,  $\alpha_t(\mathbf 1)=\sum_n\phi(e_n(t))\phi(e_n(t))^*=\mathbf 1$, 
$t> 0$.   
Thus we may conclude that the given product system 
$E$ can be associated with an $E_0$-semigroup.

\section{Proof of Theorem \ref{esThm1}}\label{S:re}

The following observation implies  that we could just as well have defined the inner product 
of (\ref{ssEq2}) by 
$$
\langle f,g\rangle=\lim_{T\to\infty}\int_T^{T+1}\langle f(\lambda),g(\lambda)\rangle\,d\lambda.
$$

\begin{lem}\label{esLem0}
For any two stable sections $f,g$, there is a $\lambda_0>0$ such that 
$$
\langle f,g\rangle = \int_T^{T+1} \langle f(\lambda),g(\lambda)\rangle\,d\lambda
$$
for all real numbers $T\geq \lambda_0$.  
\end{lem}

\begin{proof}  Let $u:(0,\infty)\to\mathbb C$ 
be a Borel function satisfying $\int_0^T|u(\lambda)|\,d\lambda<\infty$ 
for every $T>0$, together with $u(\lambda+1)=u(\lambda)$ a.e. for sufficiently 
large $\lambda$.  Then 
for $k\in\mathbb N$, the integral 
$
\int_k^{k+1}u(\lambda)\,d\lambda
$
becomes independent of $k$ when $k$ is large.  
We claim that for sufficiently 
large $T$ and the integer $n=n_T$ satisfying $T<n\leq T+1$, one has 
\begin{equation}\label{ssEq3}
\int_T^{T+1}u(\lambda)\,d\lambda = 
\int_n^{n+1}u(\lambda)\,d\lambda.  
\end{equation}
Note that Lemma \ref{esLem0} follows from (\ref{ssEq3}) after 
taking $u(\lambda)=\langle f(\lambda),g(\lambda)\rangle$.  

Of course, the formula (\ref{ssEq3}) is completely elementary. 
The integral on the left decomposes into a sum $\int_T^n \ +\  \int_n^{T+1}$, 
and for large $T$ we can write 
$$
\int_T^{n}u(\lambda)\,d\lambda=
\int_{T}^nu(\lambda+1)\,d\lambda=
\int_{T+1}^{n+1}u(\lambda)\,d\lambda. 
$$
It follows that 
\begin{align*}
\int_T^{T+1}u(\lambda)\,d\lambda=
(  \int_{T+1}^{n+1} \ +\  \int_n^{T+1})\,u(\lambda)\,d\lambda
=\int_{n}^{n+1}u(\lambda)\,d\lambda, 
\end{align*}
which proves (\ref{ssEq3}).  
\end{proof}

To show that $\phi$ is a representation, we must show that for every 
$t>0$, every $v,w\in E(t)$, and every $f,g\in \mathcal S$ one has 
$\langle\phi_0(v)f,\phi_0(w)g\rangle=\langle v,w\rangle\langle f,g\rangle$.   
Indeed, for sufficiently large $n\in\mathbb N$ we can write  
\begin{align*}
\langle \phi_0(v)f,\phi_0(w)g\rangle &=
\int_n^{n+1}\langle \phi_0(v)f(\lambda),\phi_0(w)g(\lambda)\rangle\,d\lambda \\
&=
\int_n^{n+1}\langle v\cdot f(\lambda-t),w\cdot g(\lambda-t)\rangle \,d\lambda \\
&=\langle v,w\rangle\int_n^{n+1}\langle f(\lambda-t),g(\lambda-t)\rangle\,d\lambda \\
&=\langle v,w\rangle \int_{n-t}^{n-t+1}\langle f(\lambda),g(\lambda)\rangle\,d\lambda 
=\langle v,w\rangle \langle f,g\rangle,  
\end{align*}
where the final equality uses Lemma \ref{esLem0}.  

It remains to show that $\phi$ is an essential representation, and for that, we must 
calculate the adjoints of operators in $\phi(E)$.  The following notation 
from \cite{arvMono} will be convenient.

\begin{rem}
Fix $s>0$ and 
an element $v\in E(s)$; for every $t>0$ we consider the left multiplication operator 
$\ell_v:x\in E(t)\mapsto v\cdot x\in E(s+t)$.  This operator 
has an adjoint $\ell_v^*: E(s+t)\to E(s)$, 
which we write more simply as $v^*\eta=\ell_v^*\eta$, $\eta\in E(s+t)$.  
Equivalently, for $s<t$, $v\in E(s)$, $y\in E(t)$, we write $v^*y$ for 
$\ell_v^*y\in E(t-s)$.  Note that $v^*y$ is undefined for $v\in E(s)$ and $y\in E(t)$ 
when $t\leq s$.  
\end{rem}

Given elements $u\in E(r)$, $v\in E(s)$, $w\in E(t)$, the ``associative law" 
\begin{equation}\label{esEq1}
u^*(v\cdot w)=(u^*v)\cdot w
\end{equation}
makes sense when $r\leq s$ ($t>0$ can be arbitrary), provided 
that it is suitably interpreted when $r=s$.  
Indeed, it is true {\em verbatim} when $r<s$ and $t>0$, while 
if $s=r$ and $t>0$, then it 
takes the form 
\begin{equation}\label{esEq2}
u^*(v\cdot w)=\langle v,u\rangle_{E(s)} \cdot w,\qquad u,v\in E(s),\quad w\in E(t).  
\end{equation}

\begin{lem}\label{esLem1}  Choose $v\in E(t)$. 
For every stable section $f\in\mathcal S$, there is a null section $g\in\mathcal N$ 
such that 
$$
(\phi_0(v)^*f)(\lambda)=v^*f(\lambda+t) + g(\lambda), \qquad \lambda>0.  
$$
\end{lem}

\begin{proof}
A straightforward calculation of the adjoint of 
$\phi_0(v):\mathcal S\to\mathcal S$ with respect to the semidefinite 
inner product (\ref{ssEq2}). 
\end{proof}

Lemma \ref{esLem2} follows from the identification 
$E(t)\cong E(s)\otimes E(t-s)$ when $s<t$.  We include a proof  
for completeness.  

\begin{lem}\label{esLem2}
Let $0<s<t$, let 
$v_1, v_2, \dots$ be an orthornormal basis for $E(s)$ and let $\xi\in E(t)$.  
Then 
\begin{equation}\label{esEq3.5}
\sum_{n=1}^\infty \|v_n^*\xi\|^2=\|\xi\|^2.  
\end{equation}
\end{lem}

\begin{proof}
For $n\geq 1$, $\xi\in E(t)\mapsto v_n(v_n^*\xi)\in E(t)$ 
defines a sequence of mutually orthogonal projections in $\mathcal B(E(t))$.  
We claim that these projections sum to 
the identity.  
Indeed, since $E(t)$ is the closed linear span of the set 
of products $E(s)E(t-s)$, it suffices to show that for every vector in $E(t)$ of the form 
$\xi = \eta\cdot \zeta$ with $\eta\in E(s)$, $\zeta\in E(t-s)$, we have 
$\sum_n v_n(v_n^*\xi)=\xi$.  For that, we can use 
(\ref{esEq1}) and (\ref{esEq2}) to write 
$$
v_n(v_n^*\xi)= v_n(v_n^*(\eta\cdot \zeta))=v_n((v_n^*\eta)\cdot \zeta)=\langle \eta,v_n\rangle v_n\cdot \zeta, 
$$
hence 
$$
\sum_{n=1}^\infty v_n(v_n^*\xi)=(\sum_{n=1}^\infty \langle \eta,v_n\rangle v_n)\cdot \zeta
=\eta\cdot \zeta=\xi, 
$$
as asserted.  (\ref{esEq3.5}) follows after taking the inner product with $\xi$.  
\end{proof}

\begin{proof}[Proof of Theorem \ref{esThm1}]
Since the subspaces $H_t=[\phi(E(t))H]$ satisfy 
$H_{s+t}=[\phi(E(s))H_t]\subseteq H_t$, 
it suffices to show that $H_1=H$.  For that, it is enough to show that for $\xi\in H$ of the form 
$\xi=f+\mathcal N$ where $f$ is a stable section
\begin{equation}\label{esEq3}
\langle\sum_{n=1}^\infty\phi(v_n)\phi(v_n)^*\xi,\xi\rangle =\sum_{n=1}^\infty \|\phi_0(v_n)^*f\|^2 
=\|f\|^2=\|\xi\|^2,  
\end{equation}
$v_1, v_2,\dots$ denoting an orthonormal basis for $E(1)$.  Fix such a basis 
$(v_n)$ for $E(1)$ and a stable section $f$.  Choose $\lambda_0>1$ so that 
$f(\lambda+1)=f(\lambda)\cdot e$ (a.e.) for $\lambda> \lambda_0$.  
For $\lambda>\lambda_0$ we have $\lambda+1>1$, so Lemma \ref{esLem2} implies   
$$
\sum_{n=1}^\infty \|v_n^*f(\lambda+1)\|^2=\|f(\lambda+1)\|^2=
\|f(\lambda)\cdot e\|^2=\|f(\lambda)\|^2, \quad \text{(a.e.)}.  
$$ 
It follows that for every integer $N>\lambda_0$, 
\begin{align*}
\sum_{n=1}^\infty\int_N^{N+1}\|v_n^*f(\lambda+1)\|^2\,d\lambda &=
\int_N^{N+1}\sum_{n=1}^\infty\|v_n^*f(\lambda+1)\|^2\,d\lambda 
\\&=
\int_N^{N+1}\|f(\lambda)\|^2\,d\lambda =\|f+\mathcal N\|_H^2.  
\end{align*}
Lemma \ref{esLem1} implies that when $N$ is sufficiently large, the left 
side is 
$$
\sum_{n=1}^\infty \int_N^{N+1}\|(\phi_0(v_n)^*f)(\lambda)\|^2\,d\lambda =
\sum_{n=1}^\infty \|\phi_0(v_n)f\|^2, 
$$
and (\ref{esEq3}) follows.  
\end{proof}

\begin{rem}[Nontriviality of $H$]\label{rem1}
Let $L^2((0,1];E)$ be the subspace of $L^2(E)$ consisting of all 
sections that vanish almost everywhere outside the unit interval.  Every 
$f\in L^2((0,1];E)$ corresponds to a stable section $\tilde f\in\mathcal S$ by 
extending it from $(0,1]$ to $(0,\infty)$ by periodicity
$$
\tilde f(\lambda)=f(\lambda-n)\cdot e^n,\qquad n<\lambda\leq n+1, \quad n=1,2,\dots, 
$$
and for every $n=1,2,\dots$ we have 
$$
\int_{n}^{n+1}\|\tilde f(\lambda)\|^2\,d\lambda=\int_n^{n+1}\|f(\lambda-n)\cdot e^n\|^2\,d\lambda
=\int_0^1\|f(\lambda)\|^2\,d\lambda.  
$$
Hence the map $f\mapsto \tilde f+\mathcal N$ embeds $L^2((0,1];E)$ isometrically 
as a subspace of $H$; in particular, 
$H$ is not the trivial Hilbert space $ \{0\}$.  
\end{rem}

\begin{rem}[Purity]
An $E_0$-semigroup $\alpha=\{\alpha_t: t\geq 0\}$ is said to be {\em pure} if 
the decreasing von Neumann algebras $\alpha_t(\mathcal B(H))$ 
have trivial intersection $\mathbb C\cdot \mathbf 1$.  
The question of whether every $E_0$-semigroup is a cocycle 
perturbation of a pure one has been resistant \cite{arvMono}.  
Equivalently, is every product system associated with a {\em pure} $E_0$-semigroup?  
While the answer is yes for product systems of type $I$ and $II$, and it is yes 
for the type $III$ examples constructed by Powers (see 
\cite{powTypeIII} or Chapter 13 of \cite{arvMono}), 
it is unknown in general.  

It is perhaps worth pointing out that we have shown that 
the examples of Theorem \ref{esThm1} are 
not pure; hence the above construction appears to be inadequate for 
approaching that issue.  Since the proof establishes a negative 
result that is 
peripheral to the direction of this note, 
we have omitted it.    
\end{rem}

\bibliographystyle{alpha}

\newcommand{\noopsort}[1]{} \newcommand{\printfirst}[2]{#1}
  \newcommand{\singleletter}[1]{#1} \newcommand{\switchargs}[2]{#2#1}

\end{document}